\documentclass[
submission
]{dmtcs-episciences}

\usepackage[utf8]{inputenc}
\usepackage{subfigure}

\def\0{\emptyset}

\def\n{\noindent}

\usepackage{graphicx}
 \usepackage{amsmath}
 \usepackage{amsfonts}
\usepackage{amssymb}

 \usepackage{algorithm}
 \usepackage{algorithmic}
\author{Ke Liu\thanks{liuke17@mails.tsinghua.edu.cn}
  \and Mei Lu\thanks{lumei@tsinghua.edu.cn}}
\title[The treewidth of 2-section of hypergraphs]{The treewidth of 2-section of hypergraphs}
\affiliation{
  Department of Mathematical Sciences, Tsinghua
University, Beijing 100084, China.}
\received{2020-05-24}
\revised{2021-06-16}
\accepted{2021-06-30}
\begin{document}
\publicationdetails{23}{2021}{3}{1}{6499}
\maketitle
\begin{abstract}
 Let $H=(V,F)$ be a simple hypergraph without  loops.  $H$ is called  linear if $|f\cap g|\le 1$ for any $f,g\in F$ with $f\not=g$. The $2$-section of  $H$, denoted by $[H]_2$, is a graph with $V([H]_2)=V$ and for any $ u,v\in V([H]_2)$, $uv\in E([H]_2)$ if and only if there is $ f\in F$ such that $u,v\in f$. The treewidth of a graph is an important invariant in structural and algorithmic graph theory. In this paper, we consider the treewidth of the $2$-section  of a linear hypergraph. We will use the minimum degree, maximum degree, anti-rank  and average rank of a linear hypergraph to determine the upper and lower bounds of the treewidth of its $2$-section. Since for any graph $G$, there is a linear hypergraph $H$ such that $[H]_2\cong G$, we provide a method to estimate the  bound of treewidth of graph by  the parameters of the hypergraph.

\end{abstract}

\section{Introduction}

The treewidth of a graph is an important invariant in structural and algorithmic graph theory. The concept of treewidth was originally introduced by Bertel\'e and  Brioschi \cite{BB} under the name of dimension. It was later rediscovered by Halin \cite{Hal} in 1976 and by  Robertson and  Seymour \cite{RS84} in 1984. Now it has  been studied by many other authors (see for example \cite{GO}-\cite{KS}). Treewidth is commonly used as a parameter in the parameterized complexity analysis of graph algorithms, since many NP-complete problems can be solved in polynomial time on graphs of bounded treewidth \cite{Bo}. The relation between the treewidth and other graph parameters has been explored in a number of papers (see \cite{HW17} for a recent survey).
In \cite{Har}, Harvey and Wood studied the treewidth of line graphs. They  proved sharp lower bounds of the treewidth of the line graph of a graph $G$ in terms of both the minimum degree and the average degree of $G$. Motivated by their work, in this paper, we study the treewidth of $2$-section of linear hypergraphs. 

A {\em hypergraph} is a pair $H = (V,F)$, where $V $ is a finite set
of vertices and $F$ is a family of subsets of $V$ such that for any $f\in F$, $f\neq \emptyset$ and $V=\cup_{f\in F}f$. The size of $H$ is the cardinality of $F$. A {\em simple hypergraph} is a hypergraph $H$ such that
if $f\subseteq g$, then $f=g$, where $f,g\in F$. If $|f|=1$, we call $f$ a loop. In this paper, we just consider  simple and  no loop hypergraphs.
The {\em rank} and   {\em anti-rank} of $H$ is defined as $r(H)=\max\limits_{f\in F}|f|$ and $s(H)=\min\limits_{f\in F}|f|$, respectively. If $r(H)=s(H)=2$, then $H$ is a graph. For any $v\in V$, denote $F(v)=\{f\in F|v\in f\}$. Then the degree of $v$, denoted by $deg(v)$, is  $|F(v)|$. The maximum and minimum degree of $H$ will be denoted by $\Delta=\max_{v\in V} deg(v)$ and $\delta=\min_{v\in V} deg(v)$, respectively. If $\delta=\Delta=k$, then we call the hypergraph {\em $k$-regular}. The {\em average rank} of $H$ is defined as $l (H)=\sum_{f\in F}|f|/|F|$.
A hypergraph $H$ is called {\em linear} if $|f\cap g|\le 1$ for any $f,g\in F$ with $f\not=g$.

Let $H=(V,F)$ be a hypergraph (or graph), where $V=\{v_1,\ldots,v_n\}$ and $F=\{f_1,f_2,\ldots,f_m\}$. The {\em dual } of  $H$, denoted by $H^*=(V^*,F^*)$, is a hypergraph whose vertices  $u_1,u_2,\ldots,u_m$ correspond to the edges of $H$ and with edges $g_i=\{u_j|v_i\in f_j\}$, $1\le i\le n$.
The {\em line graph} of a (hyper)graph  $H$, denoted by $L(H)$ is a graph whose vertices $w_1,w_2,\ldots,w_m$ correspond to the edges of $H$ and with edges $w_{i}w_{j}$ if $f_{i}\cap f_{j}\neq \emptyset$.
Terminology and notation concerning hypergraph not defined here can be found in \cite{Berge}. Now we give the definitions of treewidth.
Let $T$ be a tree. We will use $T$ to denote the vertex set of $T$ for short.
\vskip.2cm

{\noindent\bf Definition 1.1 }
A tree decomposition of a graph $G=(V,E)$ is a pair $(T,(B_{t})_{t\in T})$, where $T$ is a tree and $(B_{t})_{t\in T}$ ($B_t$ is called a bag) is a family of subsets of $V$ such that:

\n (T1) for every $v\in V$, the set $B^{-1}(v)=\{t\in T|v\in B_{t}\}$ is nonempty and connected in $T$;

\n (T2) for every edge $uw\in E(G)$, there is $t\in T$ such that $u,w\in B_{t}$.
\vskip.2cm

\n The width of the decomposition $(T,(B_{t})_{t\in T})$ is the number
$$\max\{\left|B_{t}\right||t\in T\}-1.$$
The {\em treewidth} $tw(G)$ of $G$ is the minimum of the widths of the tree decompositions of $G$.

\vskip.2cm
In order to  cite the definition of a generalized hypertree decomposition of a hypergraph which was given in \cite{GNF}, we introduce the definition of the 2-section of a hypergraph.
\vskip.2cm

{\noindent\bf Definition 1.2 }
   The 2-section of a hypergraph $H=(V,F)$, denoted by $[H]_2$, is a graph with $V([H]_2)=V$ and for any $ u,v\in V([H]_2)$, $uv\in E([H]_2)$ if and only if there is $ f\in F$ such that $u,v\in f$.
\vskip.2cm

{\noindent\bf Definition 1.3 \cite{GNF}}
A generalized hypertree decomposition of a hypergraph $H=(V,F)$ is a 3-tuple $(T,(B_{t})_{t\in T},(\lambda_{t})_{t\in T})$, where $T$ is a tree, $(B_{t})_{t\in T}$ is a family of subsets of $V$ and $(\lambda_{t})_{t\in T}$ is a family of subsets of $F$ such that:

\n (TI) $(T,(B_{t})_{t\in T})$ is a tree decomposition of $[H]_2$;

\n (TII) for every $t\in T$, $B_{t}\subseteq \cup_{f\in \lambda_{t}}f$.

\vskip.2cm

Now we give a new definition of decomposition of a hypergraph which is called a supertree decomposition of a hypergraph.

{\noindent\bf Definition 1.4 } A supertree decomposition of a hypergraph $H=(V,F)$ is a 3-tuple $(T,(B_{t})_{t\in T},(\lambda_{t})_{t\in T})$, where $(T,(B_{t})_{t\in T},(\lambda_{t})_{t\in T})$ is a generalized hypertree decomposition of  $H$ such that:

\n (TIII) for every $f\in F$, the set $\lambda^{-1}(f)=\{t\in T|f\in \lambda_{t}\}$ is nonempty and connected in $T$;

\n (TIV) for every $f_{1},f_{2}\in F$ with $f_{1}\cap f_{2}\neq \emptyset$, there is a $t\in T$ such that $f_{1},f_{2}\in \lambda_{t}$.

\n The width of the decomposition $(T,(B_{t})_{t\in T},(\lambda_{t})_{t\in T})$ of $H$ is the number
$$\max\{\left|\lambda_{t}\right||t\in T\}.$$
The {\em supertree width} $stw(H)$ of $H$ is the minimum of the widths of the supertree decompositions of $H$.

By Definition 1.4, if $(T,(B_{t})_{t\in T},(\lambda_{t})_{t\in T})$ is a supertree decomposition of $H$, then $(T,(\lambda_{t})_{t\in T})$ is a tree decomposition of $L(H)$. When $(T,(\lambda_{t})_{t\in T})$ is a tree decomposition of $L(H)$, we let $B_{t}=\cup_{f\in \lambda_{t}}f$. Then $(T,(B_{t})_{t\in T},(\lambda_{t})_{t\in T})$ is a supertree decomposition of $H$. So we have $stw(H)=tw(L(H))+1$.
\vskip.2cm

From the Definition 1.2, we can get the following lemma.

{\noindent\bf Lemma 1.1 }{\em Let $H$ be a $2$-regular linear hypergraph. Then $[H]_2\cong L(H^*)$.}

\begin{proof} Let $H=(V,F)$ be a $2$-regular linear hypergraph. Then $H^*$ is a simple graph. By the definition of dual hypergraph, there is a bijection $\sigma$ between the edge set $F(H)$ (resp. the vertex set $V(H)$) and the vertex set $V(H^*)$ (resp. the edge set of $E(H^*)$) such that for any $f,g\in F(H)$ (resp. for any $u,v\in V(H)$), $\sigma(f)\sigma(g)\in E(H^*)$ if and only if $f\cap g\not=\emptyset$ (resp. $\sigma(u)\cap \sigma(v)\not=\emptyset$ if and only if there is $f\in F(H)$ with $u,v\in f$). By the definition of line graph, there is a bijection $\phi:E(H^*)\rightarrow V(L(H^*))$ such that for any $e_1,e_2\in E(H^*)$, $\phi(e_1)\phi(e_2)\in E(L(H^*))$ if and only if $e_1\cap e_2\not=\emptyset$.

We will show that $[H]_2\cong L(H^*)$. Let $\phi\sigma: V([H]_2)\rightarrow V(L(H^*))$. Then $\phi\sigma$ is a bijection. For any $u,v\in V([H]_2)$, $uv\in E([H]_2)$ if and only if there is $f\in F(H)$ such that $u,v\in f$  if and only if $\sigma(u)\cap \sigma(v)\not=\emptyset$ if only if $\phi(\sigma(u))\phi(\sigma(v))\in E(L(H^*))$. Thus $[H]_2\cong L(H^*)$.\end{proof}


\vskip.2cm
By the definitions of line graph and the dual, we can easily get the following lemma.

{\noindent\bf Lemma 1.2 }{\em Let $H$ be a $2$-regular linear hypergraph. Then $H^*\cong L(H)$.}

\vskip.2cm

Let $H$ be a  hypergraph. There are two elementary lower bounds on treewidth of $[H]_2$. First,
\begin{equation}
\begin{split}
& tw([H]_2)\geq r(H)-1
\end{split}
\end{equation}
 since the vertices in a hyperedge form a clique in $[H]_2$. Second, given a minimum width tree decomposition of $[H]_2$, replace each vertex with all the hyperedges containing the vertex to obtain a supertree decomposition of $H$. It follows that
\begin{equation}
\begin{split}
& tw([H]_2)\geq \frac{1}{\Delta}stw(H)-1.
\end{split}
\end{equation}

We prove the following lower bound on $tw([H]_2)$ in terms of the minimum degree $\delta$, maximum degree $\Delta$ and average rank $l(H)$ of a linear hypergraph $H$.
\vskip.2cm

{\noindent\bf Theorem 1.1 }
{\em Let $H$ be a linear hypergraph with minimum degree $\delta$, maximum degree $\Delta$ and average rank $l(H)$. Let
$\Delta\geq\delta\geq 2$.  Suppose $\Delta\le2\delta^{2}-2\delta$. Then
$$tw([H]_2)>\left\{
\begin{array}{ll}
\frac{(2\delta^{2}-2\delta-\Delta)l(H)^{2}+(2\Delta+4\delta^{2}-4\delta)l(H)}{4\Delta\delta(\delta-1)}-1 &\mbox{if $\delta^{2}-\delta\leq\Delta-2\Delta /l(H)$,}\\
\frac{(2\delta^{2}-2\delta-\Delta)l(H)^{2}+6\Delta l(H)-8\Delta}{4\Delta\delta(\delta-1)}-1 & \mbox{otherwise.}
\end{array}
\right.$$}


\vskip.2cm
In \cite{Har}, Harvey and Wood showed that for every graph $G$, $tw(L(G)) >\frac{1}{8} d(G)^2 + \frac{3}{4}d(G)- 2$, where $d(G)$ is the average degree of $G$. Let $H$ be a 2-regular linear hypergraph of order $n$ and size $m$. By Lemma 1.1, $[H]_2\cong L(H^*)$. Note that $d(H^*)=\frac{2n}{m}=l(H)$. By Theorem 1.1, $tw(L(H^*))=tw([H]_2)> \frac{1}{8} d(H^*)^2 + \frac{3}{4}d(H^*)- 2$, just  as the result in \cite{Har}.

We also prove two lower bounds on $tw([H]_2)$ in terms of the anti-rank $s(H)$ based on different condition of minimum degree of the given hypergraph $H$.
\vskip.2cm

{\noindent\bf Theorem 1.2 }
{\em For every linear hypergraph $H$ with anti-rank $s(H)$ and minimum degree $\delta\geq 3$, we have
$$tw([H]_2)\geq\left\{
\begin{array}{ll}
\frac{3}{8}s(H)^{2}+\frac{3}{4}s(H)-1  &  \mbox{when} \ s(H) \mbox{ is even,}\\
\frac{3}{8}s(H)^{2}+\frac{1}{2}s(H)-\frac{7}{8}  &  \mbox{when} \ s(H) \mbox{ is odd}.
\end{array}\right.$$}

\vskip.2cm
{\noindent\bf Theorem 1.3 }
{\em For every linear hypergraph $H$ with anti-rank $s(H)$ and minimum degree $\delta=2$, we have
$$tw([H]_2)\geq\left\{
\begin{array}{ll}
\frac{1}{4}s(H)^{2}+s(H)-1  &  \mbox{when} \ s(H) \mbox{ is even,}\\
\frac{1}{4}s(H)^{2}+s(H)-\frac{5}{4}  &  \mbox{when} \ s(H) \mbox{ is odd}.
\end{array}\right.$$}

\vskip.2cm
In \cite{Har}, Harvey and Wood showed that for every graph $G$ with minimum degree $\delta(G)$,
$$tw(L(G))\geq\left\{
\begin{array}{ll}
\frac{1}{4}\delta(G)^{2}+\delta(G)-1  &  \mbox{when} \ \delta(G) \mbox{ is even,}\\
\frac{1}{4}\delta(G)^{2}+\delta(G)-\frac{5}{4}  &  \mbox{when} \ \delta(G) \mbox{ is odd}.
\end{array}\right.$$ Let $H$ be a 2-regular linear hypergraph. By Lemma 1.1, $[H]_2\cong L(H^*)$. Then $tw(L(H^*))$ $=tw([H]_2)$. Note that $\delta(H^*)=s(H)$. Thus we can obtain the same result as that in \cite{Har}  by Theorem 1.3.
\vskip.2cm
 Now we consider upper bounds on $tw([H]_2)$. It is easy to show that
\begin{equation}
\begin{split}
& tw([H]_2)\leq r(H)stw(H)-1.
\end{split}
\end{equation}

To see this, we consider a minimum width  supertree decomposition of $H$, and replace each bag $\lambda_{t}$ by the vertices that are incident to an hyperedge of $\lambda_{t}$. This creates a tree decomposition of $[H]_2$, where each bag contains at most $r(H)stw(H)$ vertices. In Section 5, we improve this bound as follows.

\vskip.2cm
{\noindent\bf Theorem 1.4 }
{\em For every linear hypergraph $H$, we have
$$tw([H]_2)\leq\frac{2}{3}stw(H)r(H)+\frac{1}{3}(stw(H)-1)^{2}+\frac{1}{3}r(H)-1.$$}
Theorem 1.4 is of primary interest when $r(H)\gg stw(H)$, in which case the upper bound is $(\frac{2}{3}stw(H)+\frac{1}{3})r(H)$. When $r(H)<stw(H)-1$, the bound in (3) is better than that in Theorem 1.4.

In \cite{Har}, Harvey and Wood showed that for every graph $G$, $tw(L(G))\leq\frac{2}{3}(tw(G)+1)\Delta(G)+\frac{1}{3}tw(G)^{2}+\frac{1}{3}\Delta(G)-1$. Let $H$ be a 2-regular linear hypergraph. Recall that $stw(H)=tw(L(H))+1$. By Lemmas 1.1 and 1.2, $[H]_2\cong L(H^*)$ and  $H^*\cong L(H)$. Then $tw(H^*)=stw(H)-1$.
Note that $\Delta(H^*)=r(H)$. By Theorem 1.4, $tw(L(H^*))=tw([H]_2)\leq\frac{2}{3}(tw(H^*)+1)\Delta(H^*)+\frac{1}{3}tw(H^*)^{2}+\frac{1}{3}\Delta(H^*)-1$, just the same as that in \cite{Har}.

The rest of this paper is organized as follows. In Section 2, some properties of tree decompositions of 2-section of hypergraphs are given. The proof of Theorem 1.1 is given in Section 3. In Section 4, we will prove Theorems 1.2 and  1.3. In Section 5, we prove Theorem 1.4. Section 6 concludes this paper.

\section{ Tree decomposition of 2-section of hypergraphs}

In this section, we first give  some properties of tree decompositions of 2-section of hypergraphs which will be used in next sections.

Let  $H=(V,F)$ be a hypergraph. For any $v\in V$, recall that $F(v)=\{f\in F|v\in f\}$.  Let $(T,(B_{t})_{t\in T})$ be a  tree decomposition of $[H]_2$. For $u,v\in T$, we use  $Path(u,v)$ to denote the path in $T$ connecting $u$ and $v$.

{\noindent\bf Lemma 2.1 }
{\em For every hypergraph $H=(V,F)$, there exists a minimum width tree decomposition $(T,(B_{t})_{t\in T})$ of $[H]_2$ together with an assignment $b:F\rightarrow T$ such that for each vertex $v\in V$, $B^{-1}(v)=V(ST_v)$, where $ST_v$ is the subtree of $T$ induced by $\cup_{f_i,f_j\in F(v)}Path(b(f_{i}),b(f_{j}))$.}

\begin{proof} Let $(T,(B_{t})_{t\in T})$ be a minimum width tree decomposition of $[H]_2$ such that
$\sum_{v\in V}\left|B^{-1}(v)\right|$
 is minimized. For each  $f\in F$, the vertices in $f$ form a clique in $[H]_2$. Thus there exists a bag  $B_t$ containing all the vertices in $f$, where $t\in T$. Hence for each $f\in F$ choose one such node and declare it $b(f)$.

Let $v\in V$. For any $f\in F(v)$, we have $v\in B_{b(f)}$. It follows that $V(ST_v)\subseteq B^{-1}(v)$. If $|V(ST_v)|<|B^{-1}(v)|$, then we remove $v$ from all bags $B_t$ for $t\in B^{-1}(v)\setminus V(ST_v)$. Since each vertex incident to $v$ appears in $\cup_{f\in F(v)}B_{b(f)}$, the removal of the aforementioned vertices yields another tree decomposition of $[H]_2$. However, the existence of such a tree decomposition would contradict our choice of $(T,(B_{t})_{t\in T})$. Hence $V(ST_v)=B^{-1}(v)$, as required.\end{proof}

\vskip.2cm

We call $b(f)$ the {\em base node} of $f$. From the proof of Lemma 2.1, we can construct a tree decomposition  of $[H]_2$ so that we can assign a base node for each $f\in F$ and all the vertices in $f$ are placed in the corresponding bag $B_{b(f)}$. In fact, we can obtain a slightly stronger result that will be used to prove our  theorems.

Given $(T,(B_{t})_{t\in T})$ and $b$ guaranteed by Lemma 2.1, we can also ensure that each base node is a leaf and that $b$ is a bijection between edges of $H$ and leaves of $T$. If $b(f)$ is not a leaf, then we add a leaf adjacent to $b(f)$ and let $b(f)$ be this leaf instead. If some leaf $t$ is the base node for several edges of $H$, then we add a leaf adjacent to $t$ for each edge assigned to $t$. Finally, if $t$ is a leaf that is not a base node, then delete $t$; this maintains the desired properties since a leaf is never an internal node of a subtree.

We can improve this further. Given a tree $T$, we can root it at a node and orient all edges away from the root. In such a tree, a leaf is a node with outdegree 0. Say a rooted tree is binary if every non-leaf node has outdegree 2.
Given a tree decomposition, by the same way described in \cite{Har}, we can root it and then modify the underlying tree so that each non-leaf node has outdegree 2. The above results give the following lemma.

{\noindent\bf Lemma 2.2 }
{\em For every hypergraph $H=(V,F)$ there exists a minimum width tree decomposition $(T,(B_{t})_{t\in T})$ of $[H]_2$ together with an assignment $b:F\rightarrow T$ such that
$T$ is a binary tree, $b$ is a bijection onto the leaves of $T$ and for each $v\in V$, $B^{-1}(v)=V(ST_v)$.}

\vskip.2cm
By Lemma 2.2, we have the following lower bound on $tw([H]_2)$ that is slightly stronger than (2).

\vskip.2cm
{\noindent\bf Lemma 2.3} {\em Let $H=(V,F)$ be a hypergraph with  $\Delta\ge 2$. Then we have
$tw([H]_2)\geq \frac{1}{\Delta -1}(stw(H)-1)-1$.}

\begin{proof} Let $k=tw([H]_2)+1$ and  $(T,(B_{t})_{t\in T})$  be a tree decomposition of $[H]_2$ of width $k-1$, together with an assignment $b$ as ensured by Lemma 2.2. For each $v\in V$ with $deg(v)\ge 2$, we can assume that $|B^{-1}(v)|\ge 2$; otherwise, we can simply add a leaf $t'$ adjacent to $t$ and let $B_{t'}=B_{t}$, where $t\in B^{-1}(v)$.
We now partially construct a  supertree decomposition $(T',(B'_{t})_{t\in T'},(\lambda_{t})_{t\in T'})$ of $H$ as follows: first let $(T',(B'_{t})_{t\in T'})=(T,(B_{t})_{t\in T})$.
Let $v\in V$. If $deg(v)=1$, then we just put the hyperedge adjacent to $v$ in $\lambda_t$, where $t\in B^{-1}(v)$. Assume that $deg(v)\ge 2$. We arbitrarily choose two hyperedges in $F(v)$, say $f_{v}^1$ and $f_{v}^2$. We put $F(v)\setminus  \{f_{v}^2\}$ in $\lambda_t$ for all $t\in B^{-1}(v)\setminus \{b(f_{v}^2)\}$ and put  $F(v)\setminus \{f_{v}^1\}$ in $\lambda_{b(f_{v}^2)}$. Then for all $t$, $|\lambda_{t}|\le k(\Delta-1)$ since each vertex contributes at most $\Delta-1$ hyperedges to a given bag. An edge $tt'$ in $T$ is called the {\em edge corresponding to} a 3-tuple $(v,f_{v}^1,f_{v}^2)$ if $f_{v}^1\in \lambda_{t}\setminus \lambda_{t'}$ and $f_{v}^2\in \lambda_{t'}\setminus \lambda_{t}$.
If $tt'$ is an edge corresponding to  $(v,f_{v}^1,f_{v}^2)$ and $(v',f_{v'}^1,f_{v'}^2)$ simultaneously, say $f_{v}^1,f_{v'}^1\in \lambda_{t}$, then we subdivide $tt'$ by adding
 a new node $t''$ and let $B'_{t''}=B'_{t}\cap B'_{t'},\lambda_{t''}=(\lambda_t\setminus \{f_{v}^1\})\cup\{f_{v}^2\}$. Repeat this process so that every edge in $T'$ corresponds to at most one 3-tuple. We also use $T'$ to denote the tree after all of these subdivisions.
 Note that  $(T',(B'_{t})_{t\in T'})$ is still a tree decomposition of $[H]_2$.  Now if there is an edge $tt'$ corresponding to $(v,f_{v}^1,f_{v}^2)$, then we add $f_{v}^1$ into $\lambda_{t'}$, which increases the size of  $\lambda'_{t}$  at most 1.

We are going to show that $(T',(B'_{t})_{t\in T'},(\lambda_{t})_{t\in T'})$ is a supertree decomposition of $H$. By the construction of $(T',(B'_{t})_{t\in T'},(\lambda_{t})_{t\in T'})$, (TI) and (TIV) in Definitions 2.3 and 2.4 hold. So we just need to show (TII) and (TIII).

First, we prove that for all $t\in T'$, $B'_{t}\subseteq \cup_{f\in \lambda_{t}}f$. If $t\in T$, then the conclusion holds because for each $v\in B_{t}$, we put at least one edge incident to $v$ in $\lambda_{t}$. If $t$ is the subdivided node, then the result holds by the construction.

Now we show that for any $f\in F$, $\lambda^{-1}(f)$ is connected in $T'$. Suppose there is $f\in F$ such that $\lambda^{-1}(f)$ is nonconnected in $T'$. Then there must be $t,t'\in T'$ such that for all  $ t''\in Path(t,t')\setminus \{t,t'\}$, $f\in (\lambda_{t}\cap \lambda_{t'})\setminus \lambda_{t''}$. If there is a vertex $v\in f$ such that $Path(t,t')\subseteq B^{-1}(v)$, then we should have $f\in \lambda_{t''}$ for all $ t''\in Path(t,t')$ even if $t''$ is a subdivided node by the construction, a contradiction. Suppose  there exist two vertices $v,v'\in f$ such that $t\in B^{-1}(v)\setminus B^{-1}(v'),t'\in B^{-1}(v')\setminus B^{-1}(v)$, then there is $ t''\in Path(t,t')\setminus\{t,t'\}$ such that $t''\in B^{-1}(v)\cap B^{-1}(v')$ by $vv'\in E([H]_2)$ which implies that $f\in \lambda_{t''}$,  a contradiction.

Thus we have $stw(H)\leq (\Delta-1) k +1=(\Delta-1)(tw([H]_2)+1)+1$, as required. \end{proof}

\section{ Lower bound in terms of average rank}

In this section, we prove Theorem 1.1. Let $H=(V,F)$ be a hypergraph with size $m$. Let $V_{i}=\{v\in V(H)|deg(v)=i\}$ and $n_i=|V_i|$, where $i=\delta,\ldots,\Delta$. By the definition of average rank, $l(H)=(\delta n_{\delta}+\ldots+\Delta n_{\Delta})/m$.
Given two sets $X,Y\subseteq F$. Let $\sigma^{j}_{i}(X)=\#\{v\in V_i|deg_X(v)=j\}$ and let $\sigma_{i}^{j,l}(X,Y)=\#\{v\in V_i|deg_X(v)=j\mbox{~and~}deg_Y(v)=l\}$, where $deg_X(v)=\#\{f|f\in F(v)\cap X\}$. We say a hypergraph $H$ is {\em minimal} if $l(H_{S})<l(H)$ for every nonempty proper subset $S$ of $F$, where $$l(H_S)=\left[\delta\left(n_{\delta}-\sum\limits_{j=1}^{\delta}\sigma^{j}_{\delta}(S)\right)+\ldots+\Delta\left(n_{\Delta}-\sum\limits_{j=1}^{\Delta}\sigma^{j}_{\Delta}(S)\right)\right]/(m-|S|).$$
Firstly, we need some lemmas before we start bounding the tree width of $[H]_2$.

\vskip.2cm
{\noindent\bf Lemma 3.1 }
{\em If $H$ is a minimal hypergraph and $S$ is a nonempty proper subset of $F(H)$, then
$$\frac{1}{\Delta}l(H)<\frac{1}{|S|}\left(\sum\limits_{f\in S}|f|-
\sum\limits_{i=\delta}^{\Delta}\sum\limits_{j=1}^{i}\sigma^{j}_{i}(S)\left(j-\frac{i}{\Delta}\right)\right).$$}
\begin{proof} Since $H$ is minimal, we have $l(H)>l(H_{S})$ which implies
$$\frac{\sum\limits_{i=\delta}^{\Delta}i  n_{i}}{m}>
\frac{\sum\limits_{i=\delta}^{\Delta}i n_{i}-
\sum\limits_{i=\delta}^{\Delta}i \sum\limits_{j=1}^{i}\sigma^{j}_{i}(S)}{m-|S|}.$$
Hence
$l(H)=\left(\sum\limits_{i=\delta}^{\Delta}i n_{i}\right)/m<\left(\sum\limits_{i=\delta}^{\Delta}\sum\limits_{j=1}^{i}i\sigma^{j}_{i}(S)\right )/|S|.$
We  easily get $\sum\limits_{f\in S}\left|f\right|=\sum\limits_{i=\delta}^{\Delta}\sum\limits_{j=1}^{i}j\sigma^{j}_{i}(S)$. Thus
$$\frac{1}{\Delta}l(H)<
\frac{1}{|S|}\sum\limits_{i=\delta}^{\Delta}\sum\limits_{j=1}^{i}\frac{i}{\Delta}\sigma^{j}_{i}(S)
=\frac{1}{|S|}\left(\sum\limits_{f\in S}\left|f\right|-
\sum\limits_{i=\delta}^{\Delta}\sum\limits_{j=1}^{i}\sigma^{j}_{i}(S)\left(j-\frac{i}{\Delta}\right)\right),$$
and we have the conclusion.\end{proof}

\vskip.2cm
Let $(T,(B_{t})_{t\in T})$ be a tree decomposition of $[H]_2$ as guaranteed by Lemma 2.2.
For each node $t$ of $T$, let $T_{t}$ denote the subtree of $T$ rooted at $t$ containing exactly $t$ and the descendants of $t$. And let $z(T_{t})$ be the set of edges of $H$ with the base nodes in $T_{t}$. (Recall all base nodes are leaves.)

\vskip.2cm
{\noindent\bf Lemma 3.2}
{\em Let $H=(V,F)$ be a minimal hypergraph and $(T,(B_{t})_{t\in T})$ be a tree decomposition of $[H]_2$ as guaranteed by Lemma 2.2. If $t\in T$ is a non-leaf, non-root node  and $a,b$ are the children of $t$, then
$$|B_{t}|>
\frac{1}{\Delta}(|z(T_{a})|+|z(T_{b})|)l(H)-
\sum\limits_{i=\delta}^{\Delta}\sigma_{i}^{i}(z(T_{a}))-
\sum\limits_{i=\delta}^{\Delta}\sigma_{i}^{i}(z(T_{b})).$$}

\begin{proof} Denote
$$g(z(T_{a}),z(T_{b}))=\sum\limits_{f\in z(T_{a})}|f|+
\sum\limits_{f\in z(T_{b})}|f|-
\sum\limits_{i=\delta}^{\Delta}\sum\limits_{j=1}^{i}\sum\limits_{u+w=
j,u,w\geq 0}\sigma^{u,w}_{i}(z(T_{a}),z(T_{b}))\left(j-\frac{i}{\Delta}\right).$$
We first show that $g(z(T_{a}),z(T_{b}))>\frac{1}{\Delta}(|z(T_{a})|+|z(T_{b})|)l(H)$.

Since $t\in T$ is a non-leaf, non-root node,  $z(T_{a}), z(T_{b})\neq \emptyset$, $z(T_{a})\cap z(T_{b})=\emptyset$ and $z(T_{t})=z(T_{a})\cup z(T_{b})\subsetneq F$. By Lemma 3.1, we have
$$\frac{1}{\Delta}l(H)<
\frac{1}{|z(T_{a})\cup z(T_{b})|}\left(\sum\limits_{f\in z(T_{a})\cup z(T_{b})}|f|-
\sum\limits_{i=\delta}^{\Delta}\sum\limits_{j=1}^{i}\sigma^{j}_{i}(z(T_{a})\cup z(T_{b}))\left(j-\frac{i}{\Delta}\right)\right).$$
So $g(z(T_{a}),z(T_{b}))>\frac{1}{\Delta}(|z(T_{a})|+|z(T_{b})|)l(H)$.
We consider $B_{t}$, the bag of the target node $t$, which  consists of  vertices covered by edges in $z(T_{a})$ (resp. $z(T_{b})$) and $F\setminus z(T_{a})$ (resp. $F\setminus z(T_{b})$) at the same time. Thus,
\begin{eqnarray*}
|B_{t}|
&=&\sum\limits_{i=\delta}^{\Delta}\sum\limits_{j=1}^{i}\sigma^{j}_{i}(z(T_{a}))+
\sum\limits_{i=\delta}^{\Delta}\sum\limits_{j=1}^{i}\sigma^{j}_{i}(z(T_{b}))-
\sum\limits_{i=\delta}^{\Delta}\sum\limits_{j=1}^{i}\sum\limits_{u+w=j,u,w> 0}\sigma^{u,w}_{i}(z(T_{a}),z(T_{b}))\\
& &-
\sum\limits_{i=\delta}^{\Delta}\sigma_{i}^{i}(z(T_{a}))-\sum\limits_{i=\delta}^{\Delta}\sigma_{i}^{i}(z(T_{b}))\\
&=&\sum\limits_{f\in z(T_{a})}|f|-\sum\limits_{i=\delta}^{\Delta}\sum\limits_{j=1}^{i}\sigma^{j}_{i}(z(T_{a}))(j-1)+
\sum\limits_{f\in z(T_{b})}|f|-\sum\limits_{i=\delta}^{\Delta}\sum\limits_{j=1}^{i}\sigma^{j}_{i}(z(T_{b}))(j-1)\\
& &- \sum\limits_{i=\delta}^{\Delta}\sum\limits_{j=1}^{i}\sum\limits_{u+w=j,u,w> 0}\sigma^{u,w}_{i}(z(T_{a}),z(T_{b}))-
\sum\limits_{i=\delta}^{\Delta}\sigma_{i}^{i}(z(T_{a}))-\sum\limits_{i=\delta}^{\Delta}\sigma_{i}^{i}(z(T_{b})).
\end{eqnarray*}
Notice that $\sigma_{i}^{j}(z(T_{a}))=\sum\limits_{u=0}^{i-j}\sigma_{i}^{j,u}(z(T_{a}),z(T_{b}))$ and
$\sigma_{i}^{j}(z(T_{b}))=\sum\limits_{u=0}^{i-j}\sigma_{i}^{j,u}(z(T_{b}),z(T_{a}))$. We have $\sum\limits_{i=\delta}^{\Delta}\sum\limits_{j=1}^{i}\sigma_{i}^{j}(z(T_{x}))=
\sum\limits_{i=\delta}^{\Delta}\sum\limits_{j=1}^{i}\sum\limits_{u=0}^{i-j}\sigma_{i}^{j,u}(z(T_{x}),z(T_{y}))
=\sum\limits_{i=\delta}^{\Delta}\sum\limits_{j=1}^{i}\sum\limits_{u+w=j,u>0,w\ge 0}\sigma^{u,w}_{i}(z(T_{x}),z(T_{y}))$ holds for $x=a,y=b$ or $x=b,y=a$.


So the above formula can be written as
\begin{eqnarray*}
|B_{t}|&=&\sum\limits_{f\in z(T_{a})\cup z(T_{b})}|f|-
\sum\limits_{i=\delta}^{\Delta}\sum\limits_{j=1}^{i}\sum\limits_{u+w=j,u,w>0}\sigma^{u,w}_{i}(z(T_{a}),z(T_{b}))(1+u-1+w-1)\\
&-&
\sum\limits_{i=\delta}^{\Delta}\sigma_{i}^{i}(z(T_{a}))
-\sum\limits_{i=\delta}^{\Delta}\sigma_{i}^{i}(z(T_{b}))-
 \sum\limits_{i=\delta}^{\Delta}\sum\limits_{j=1}^{i}\sum\limits_{u+w=j,u\cdot w=0}\sigma^{u,w}_{i}(z(T_{a}),z(T_{b}))(j-1)\\
&=&\sum\limits_{f\in z(T_{a})\cup z(T_{b})}|f|-
\sum\limits_{i=\delta}^{\Delta}\sum\limits_{j=1}^{i}\sum\limits_{u+w=j,u,w\geq0}\sigma^{u,w}_{i}(z(T_{a}),z(T_{b}))(u-1+w-1+1)\\
&-&
\sum\limits_{i=\delta}^{\Delta}\sigma_{i}^{i}(z(T_{a}))
-\sum\limits_{i=\delta}^{\Delta}\sigma_{i}^{i}(z(T_{b}))\\
&=&\sum\limits_{f\in z(T_{a})\cup z(T_{b})}|f|-
\sum\limits_{i=\delta}^{\Delta}\sum\limits_{j=1}^{i}\sum\limits_{u+w=j,u,w\geq0}\sigma^{u,w}_{i}(z(T_{a}),z(T_{b}))(j-1)-
\sum\limits_{i=\delta}^{\Delta}\sigma_{i}^{i}(z(T_{a}))\\
&-&\sum\limits_{i=\delta}^{\Delta}\sigma_{i}^{i}(z(T_{b}))
\end{eqnarray*}
Furthermore, we have
\begin{eqnarray*}
|B_{t}|&\geq &\sum\limits_{f\in z(T_{a})\cup z(T_{b})}|f|-
\sum\limits_{i=\delta}^{\Delta}\sum\limits_{j=1}^{i}\sum\limits_{u+w=j,u,w\geq0}\sigma^{u,w}_{i}(z(T_{a}),z(T_{b}))\left(j-\frac{i}{\Delta}\right)-
\sum\limits_{i=\delta}^{\Delta}\sigma_{i}^{i}(z(T_{a}))\\
&-&\sum\limits_{i=\delta}^{\Delta}\sigma_{i}^{i}(z(T_{b}))\\
&=&g(z(T_{a}),z(T_{b}))-\sum\limits_{i=\delta}^{\Delta}\sigma_{i}^{i}(z(T_{a}))-\sum\limits_{i=\delta}^{\Delta}\sigma_{i}^{i}(z(T_{b}))\\
&> &\frac{1}{\Delta}(|z(T_{a})|+|z(T_{b})|)l(H)-
\sum\limits_{i=\delta}^{\Delta}\sigma_{i}^{i}(z(T_{a}))-
\sum\limits_{i=\delta}^{\Delta}\sigma_{i}^{i}(z(T_{b})),
\end{eqnarray*}
thus we have the conclusion. \end{proof}

\vskip.2cm
Theorem 1.1 follows from the following lemma since every hypergraph $H$ with $\delta\ge 2$ contains a minimal subgraph $H'$ with $l(H')\geq l(H)$, in which case $[H']_2\subseteq [H]_2$ and $tw([H]_2)\geq tw([H']_2)$.
\vskip.2cm
{\noindent\bf Lemma 3.3 } {\em Let $H$ be a minimal linear hypergraph with minimum degree $\delta$, maximum degree $\Delta$ and average rank $l(H)$. Let
$\Delta\geq\delta\geq 2$.  Suppose $\Delta\le2\delta^{2}-2\delta$. Then
$$tw([H]_2)>\left\{
\begin{array}{ll}
\frac{(2\delta^{2}-2\delta-\Delta)l(H)^{2}+(2\Delta+4\delta^{2}-4\delta)l(H)}{4\Delta\delta(\delta-1)}-1 &\mbox{if $\delta^{2}-\delta\leq\Delta-2\Delta /l(H)$,}\\
\frac{(2\delta^{2}-2\delta-\Delta)l(H)^{2}+6\Delta l(H)-8\Delta}{4\Delta\delta(\delta-1)}-1 & \mbox{otherwise.}
\end{array}
\right.$$}

\begin{proof} If $0\leq l(H)<2$, we have $\delta^{2}-\delta> 0\geq \Delta-2\Delta /l(H)$. Then $tw([H]_2)\geq 0\geq \frac{2}{\Delta}-1=\frac{2^{2}(2\delta^{2}-2\delta-\Delta)+6\cdot2\cdot \Delta-8\Delta}{4\Delta\delta(\delta-1)}-1>
\frac{(2\delta^{2}-2\delta-\Delta)l(H)^{2}+6\Delta l(H)-8\Delta}{4\Delta\delta(\delta-1)}-1$, as required. So we assume that $l(H)\geq 2$.  Since $\delta\ge 2$ and $H$ is a linear hypergraph, we have $l(H)<|F|$.

Let $(T,(B_{t})_{t\in T})$ be a tree decomposition of $[H]_2$ as guaranteed by Lemma 2.2. Call a node $t$ of $T$ {\em significant} if $|z(T_{t})|>\frac{l(H) }{2}$ but $|z(T_{t'})|\le\frac{l(H)}{2}$ for each child $t'$ of $t$.

\vskip.2cm
{\noindent\bf Claim 1.}
There exists a non-root, non-leaf significant node $t$.

\begin{proof} [of Claim 1] Starting at the root of $T$, begin traversing down the tree by the following rule: if some child $t$ of the current node has $|z(T_{t})|>\frac{l(H)}{2}$, then traverse to $t$; otherwise stop. Clearly this algorithm halts.

Since $|z(T_{t})|=1\le \frac{l(H)}{2}$ for any leaf $t$,  the algorithm will stop  at a non-leaf.

Let $t$ be the node where the above algorithm stops. Suppose that $t$ is the root. Then $|z(T_{t})|=|F|$. Assume $t_{1}$ and $t_{2}$ are the children of $t$. Then $|z(T_{t_{1}})|,|z(T_{t_{2}})|\leq\frac{l(H)}{2}$. Thus $|z(T_{t})|=|z(T_{t_{1}})|+|z(T_{t_{2}})|\leq l(H)<|F|$, a contradiction. Hence the algorithm does not stop at the root.\end{proof}


\vskip.2cm
Let $t$ be a non-root, non-leaf significant node and $a,b$ be the children of $t$. Set $A=z(T_{a})$ and $B=z(T_{b})$. Then $|A|,|B|\leq \frac{l(H)}{2}$ but $|z(T_{t})|=|A\cup B|> \frac{l(H)}{2}\geq1$. 
By Lemma 3.2, we can get
$$|B_{t}|> \frac{1}{\Delta}(|A|+|B|)l(H)-\sum\limits_{i=\delta}^{\Delta}\sigma_{i}^{i}(A)-
\sum\limits_{i=\delta}^{\Delta}\sigma_{i}^{i}(B).$$
We consider the following  linear programming.
\begin{equation*}
\begin{split}
 {\bf max}\quad &\sum\limits_{i=\delta}^{\Delta}\sigma_{i}^{i}(A)\\
 {\bf s.t}\quad &\sum\limits_{i=\delta}^{\Delta}\frac{i(i-1)}{2}\sigma_{i}^{i}(A)\leq \frac{|A|(|A|-1)}{2},\\
 &\sigma_{i}^{i}(A)\geq 0,\quad \delta\leq i\leq \Delta,
\end{split}
\end{equation*}
where  the constraint condition is based on $H$ being a linear hypergraph, that is, each pair ($f_{i}$, $f_{j}$) can only be calculated at most one in $\sigma_{i}^{i}(A)$. From the above linear programming, we can easily know that $\sum\limits_{i=\delta}^{\Delta}\sigma_{i}^{i}(A)\le \frac{|A|(|A|-1)}{\delta(\delta-1)}$. Thus
$$|B_{t}|>
\frac{1}{\Delta}(|A|+|B|)l(H)-
\frac{|A|(|A|-1)}{\delta(\delta-1)}-
\frac{|B|(|B|-1)}{\delta(\delta-1)}.$$
Define $\alpha,\beta,s$  such that $|A|=\alpha l(H)$, $|B|=\beta l(H)$ and $s=\frac{1}{l(H)}$. Recall $|A|,|B|\leq \frac{1}{2}l(H)$ and $|A|+|B|>\frac{1}{2}l(H)$. Hence  $|A|,|B|\geq1$. Thus $s\leq\alpha,\beta\leq\frac{1}{2}$ and $\alpha+\beta>\frac{1}{2}$. Now we have
\begin{equation*}
\begin{split}
|B_{t}|&>\frac{1}{\Delta}(\alpha l(H)+\beta l(H))l(H)-
\frac{\alpha l(H)(\alpha l(H)-1)}{\delta(\delta-1)}-
\frac{\beta l(H)(\beta l(H)-1)}{\delta(\delta-1)}\\
&=l(H)^{2}\left(\frac{\alpha}{\Delta}+\frac{\beta}{\Delta}+\frac{1}{\delta(\delta-1)}(-\alpha^{2}+s\alpha-\beta^2+s\beta)\right).\\
\end{split}
\end{equation*}
In {\bf Appendix A}, we prove that $f(\alpha,\beta)=\frac{\alpha}{\Delta}+\frac{\beta}{\Delta}+\frac{1}{\delta(\delta-1)}(-\alpha^{2}+s\alpha-\beta^2+s\beta)
\ge \min\left\{f\left(\frac{1}{2},s\right),f\left(\frac{1}{2}-s,s\right)\right\}$. Since
$f(\frac{1}{2},s)>f(\frac{1}{2}-s,s)$ if and only if
$\delta^{2}-\delta>\Delta-2\Delta s$ and $tw([H]_2)+1\geq |B_{t}|$, the result holds immediately.
\end{proof}
\vskip.2cm

By Theorem 1.1, we can get the following corollary.

{\noindent\bf Corollary 3.4 }
{\em Let $H=(V,F)$ be a $h$-regular linear hypergraph with average rank $l(H)$. Suppose $h\ge 2$.
Then
$$tw([H]_2)>\frac{(2h-3)l(H)^{2}+6l(H)-8}{4h(h-1)}-1.$$}

\vskip.2cm
Let $k\ge 2$ be a  positive integer. We construct a hypergraph $H=(V,F)$ where $F=\{f_{1},f_{2},\ldots,f_{n}\}$. The vertex set of $H$ is determined as the following rule: for any positive integers $ i,j$ with $i\not=j$, if $|i-j|\leq k$, then we put the vertex $v_{i,j}$ into both $f_{i}$ and $f_{j}$, where we let $v_{i,j}=v_{j,i}$. Then $V=\cup_{i=1}^nf_i$. We can easily know that $H$ is a $2$-regular linear hypergraph with $l(H)=2k-\gamma$, where $\gamma\to 0$ when $n\to \infty$. By Corollary 3.4, $tw([H]_2)>\frac{1}{2}k^{2}+\frac{3}{2}k-2-\gamma(\frac{1}{2}k+\frac{3}{4}-\frac{1}{8}\gamma)$. Since $\frac{1}{2}k^{2}+\frac{3}{2}p-2$ is an integer, $tw([H]_2)\geq \frac{1}{2}k^{2}+\frac{3}{2}k-2$. Note that $[H]_2\cong L(P_n^k)$, where $P_n^k$ is the $k^{th}$-power of an $n$-vertex path with vertex set $\{v_1,v_2,\ldots,v_n\}$ and edge set $\{(v_i,v_j)~|~|i-j|\leq k, 1\leq i<j\leq n\}$. In \cite{Har}, Harvey and Wood showed that $tw(L(P_n^k))\le \frac{1}{2}k^2+\frac{3}{2}k-1$. Thus when $h=2$, Corollary 3.4 is almost precisely sharp for treewidth.

\section{ Lower bounds in terms of anti-rank}

We use similar techniques to those in Section 3 to prove a lower bound on $tw([H]_2)$ in terms of $s(H)$ instead of $l(H)$.

\vskip.2cm
\begin{proof}[of Theorems 1.2 and 1.3 ]If $s(H)<2$, then Theorems 1.2 and 1.3 are trivial, since $tw([H]_2)\geq0$ whenever $[H]_2$ contains at least one vertex. Now we assume that $s(H)\geq2$. Since $H$ is a linear hypergraph and $\delta\ge 2$, $s(H)<|F|$.

Let $(T,(B_{t})_{t\in T})$ be a tree decomposition for $[H]_2$ as guaranteed by Lemma 2.2. For each node $t$ of $T$, let $T_{t}$ denote the subtree of $T$ rooted at $t$ containing exactly $t$ and the descendants of $t$. Let $z(T_{t})$ be the set of edges of $H$ with the base nodes in $T_{t}$. (Recall all base nodes are leaves.)

Call a node $t$ of $T$ significant if $|z(T_{t})|>\frac{s(H)}{2}$ but $|z(T_{t'})|\le \frac{s(H)}{2}$ for each child $t'$ of t. By a argument similar to Claim 1, there exists a non-root, non-leaf significant node $t$. Let $a,b$ be the children of $t$, and let $A=z(T_{a})$ and $B=z(T_{b})$. Then $|A|,|B|\geq 1$, $|A|,|B|\leq \frac{s(H)}{2}$ but $|A\cup B|> \frac{s(H)}{2}$.
Since $|A|,|B|$ are integers, if $s(H)$ is odd (resp. even) then $|A|+|B|\geq \frac{s(H)}{2}+\frac{1}{2}$ (resp. $|A|+|B|\geq \frac{s(H)}{2}+1$). Define $\alpha,\beta,s$ such that $|A|=\alpha s(H),|B|=\beta s(H)$ and $s=\frac{1}{s(H)}$. Then
$s\leq\alpha,\beta\leq \frac{1}{2}$ and
$$\alpha +\beta\geq\left\{
\begin{array}{ll}
\frac{1}{2}+s  &  \mbox{when} \ s(H) \mbox{ is even,}\\
\frac{1}{2}+\frac{s}{2}  &  \mbox{when} \ s(H) \mbox{ is odd}.
\end{array}\right.$$

We first give the {\bf proof of Theorem 1.2} where $\delta\ge 3$. As in Lemma 3.2,
\begin{equation*}
\begin{split}
|B_{t}|&=\sum\limits_{f\in A\cup B}|f|-
\sum\limits_{i=\delta}^{\Delta}\sum\limits_{j=1}^{i}\sum\limits_{u+w=j,u,w\geq0}\sigma^{u,w}_{i}(A,B)(j-1)-
\sum\limits_{i=\delta}^{\Delta}\sigma_{i}^{i}(A)-\sum\limits_{i=\delta}^{\Delta}\sigma_{i}^{i}(B)\\
&=\sum\limits_{f\in A\cup B}|f|-
\sum\limits_{i=\delta}^{\Delta}\sum\limits_{j=1}^{i}\sum\limits_{u+w=j,u,w\geq0}\sigma^{u,w}_{i}(A,B)(j-1)-
\sum\limits_{i=\delta}^{\Delta}\sum_{u+w=i,uw=0} \sigma_{i}^{u,w}(A,B)\end{split}
\end{equation*}
Furthermore, we have
\begin{equation*}
\begin{split}
|B_{t}|&\geq \sum\limits_{f\in A\cup B}|f|-
\sum\limits_{i=\delta}^{\Delta}\sum\limits_{j=1}^{i}\sum\limits_{u+w=j,u,w\geq0}\sigma^{u,w}_{i}(A,B)(j-1)-
\sum\limits_{i=\delta}^{\Delta}\sum_{u+w=i,u,w\geq0} \sigma_{i}^{u,w}(A,B)\\
&=\sum\limits_{f\in A\cup B}|f|-
\sum\limits_{i=\delta}^{\Delta}\sum\limits_{j=1}^{i-1}\sum\limits_{u+w=j,u,w\geq0}\sigma^{u,w}_{i}(A,B)(j-1)-
\sum\limits_{i=\delta}^{\Delta}\sum_{u+w=i,u,w\geq0} i\sigma_{i}^{u,w}(A,B)\\
&=\sum\limits_{f\in A\cup B}|f|-
\sum\limits_{i=\delta}^{\Delta}\sum\limits_{j=1}^{i-1}\sigma^{j}_{i}(A\cup B)(j-1)-
\sum\limits_{i=\delta}^{\Delta}i\sigma_{i}^{i}(A\cup B)\\
&=\sum\limits_{f\in A\cup B}|f|-
\sum\limits_{j=2}^{\Delta}\sum\limits_{i=j+1}^{\Delta}\sigma^{j}_{i}(A\cup B)(j-1)-
\sum\limits_{i=\delta}^{\Delta}i\sigma_{i}^{i}(A\cup B),
\end{split}
\end{equation*}where $\sigma^{j}_{i}(A\cup B)=0$ if $i< \delta$ or $i>\Delta$.
In {\bf Appendix B}, we show that

$$-\sum\limits_{j=2}^{\Delta}\sum\limits_{i=j+1}^{\Delta}\sigma^{j}_{i}(A\cup B)(j-1)-
\sum\limits_{i=\delta}^{\Delta}i\sigma_{i}^{i}(A\cup B)\ge -
\frac{|A\cup B|(|A\cup B|-1)}{2}.$$
Thus we have
\begin{equation*}
\begin{split}
|B_{t}|& \geq\sum\limits_{f\in A\cup B}|f|-
\frac{|A\cup B|(|A\cup B|-1)}{2}\\
& \geq(|A|+|B|)s(H)-
\frac{(|A|+|B|)(|A|+|B|-1)}{2}\\
& =\left(\left(1+\frac{1}{2}s\right)\alpha-\frac{1}{2}\alpha^{2}+\left(1+\frac{1}{2}s\right)\beta-\frac{1}{2}\beta^{2}-\alpha\beta\right)s(H)^{2}.
\end{split}
\end{equation*}
Now we calculate the minimum value of
 $f(\alpha,\beta)=
(1+\frac{1}{2}s)\alpha-\frac{1}{2}\alpha^{2}+(1+\frac{1}{2}s)\beta-\frac{1}{2}\beta^{2}-\alpha\beta$ under the condition $0\leq s\leq\alpha,\beta\leq \frac{1}{2}$ and $$\alpha +\beta\geq\left\{
\begin{array}{ll}
\frac{1}{2}+s  &  \mbox{when} \ s(H) \mbox{ is even,}\\
\frac{1}{2}+\frac{s}{2}  &  \mbox{when} \ s(H) \mbox{ is odd}.
\end{array}\right.$$Consider the partial derivative of $f(\alpha,\beta)$ with respect to $\alpha$ and $\beta$, we have
\[\frac{\partial f(\alpha,\beta)}{\partial \alpha} = \frac{s}{2}-\alpha-\beta+1>0, ~~\frac{\partial f(\alpha,\beta)}{\partial \beta} = \frac{s}{2}-\alpha-\beta+1>0.\]
Since $f(\alpha,\beta)=f(\beta,\alpha)$, we have $f(\alpha,\beta)\ge \min\left\{f\left(\frac{1}{2},\frac{1}{2}\right),f\left(\frac{1}{2},s\right)\right\}$ when $s(H)$ is even and $f(\alpha,\beta)\ge \min\left\{f\left(\frac{1}{2},\frac{1}{2}\right),f\left(\frac{1}{2},s\right),f\left(s,\frac{1}{2}-\frac{1}{2}s\right)\right\}$ when $s(H)$ is odd.
Since
$f\left(\frac{1}{2},\frac{1}{2}\right)=\frac{s}{2}+\frac{1}{2}$,
$f\left(\frac{1}{2},s\right)=\frac{3s}{4}+\frac{3}{8}$ and
$f\left(\frac{1}{2}-\frac{1}{2}s,s\right)=\frac{s^{2}}{8}+\frac{s}{2}+\frac{3}{8},$ we have
$$f(\alpha,\beta)\geq\left\{
\begin{array}{ll}
\frac{3s}{4}+\frac{3}{8}  &  \mbox{when} \ s(H) \mbox{ is even,}\\
\frac{s^{2}}{8}+\frac{s}{2}+\frac{3}{8}  &  \mbox{when} \ s(H) \mbox{ is odd}.
\end{array}\right.$$
Thus
$$tw([H]_2)+1\geq |B_{t}|\geq\left\{
\begin{array}{ll}
\frac{3}{8}s(H)^{2}+\frac{3}{4}s(H)  &  \mbox{when} \ s(H) \mbox{ is even,}\\
\frac{3}{8}s(H)^{2}+\frac{1}{2}s(H)+\frac{1}{8}  &  \mbox{when} \ s(H) \mbox{ is odd}.
\end{array}\right.$$

Now we give the {\bf Proof of Theorem 1.3} where $\delta=2$. We have
\begin{equation*}
\begin{split}
|B_{t}|&=\sum\limits_{f\in A\cup B}|f|-
\sum\limits_{i=2}^{\Delta}\sum\limits_{j=1}^{i}\sum\limits_{u+w=j,u,w\geq0}\sigma^{u,w}_{i}(A,B)(j-1)-
\sum\limits_{i=2}^{\Delta}\sigma_{i}^{i}(A)-\sum\limits_{i=2}^{\Delta}\sigma_{i}^{i}(B)\\
&=\sum\limits_{f\in A\cup B}|f|-
\sum\limits_{i=3}^{\Delta}\sum\limits_{j=1}^{i}\sum\limits_{u+w=j,u,w\geq0}\sigma^{u,w}_{i}(A,B)(j-1)-
\sum\limits_{i=3}^{\Delta}\sigma_{i}^{i}(A)-\sum\limits_{i=3}^{\Delta}\sigma_{i}^{i}(B)\\
&\quad -\sum\limits_{j=1}^{2}\sum\limits_{u+w=j,u,w\geq0}\sigma^{u,w}_{i}(A,B)(j-1)-\sigma^{2}_{2}(A)-\sigma^{2}_{2}(B)\\
&=\sum\limits_{f\in A\cup B}|f|-
\sum\limits_{j=2}^{\Delta}\sum\limits_{i=j+1}^{\Delta}(j-1)\sigma^{j}_{i}(A\cup B)-
\sum\limits_{i=3}^{\Delta}i\sigma_{i}^{i}(A\cup B)-
\sigma_{2}^{2}(A\cup B)-\sigma^{2}_{2}(A)\\
&-\sigma^{2}_{2}(B).
\end{split}
\end{equation*}
In {\bf Appendix C}, we show that
\begin{equation*}
\begin{split}&-\sum\limits_{j=2}^{\Delta}\sum\limits_{i=j+1}^{\Delta}(j-1)\sigma^{j}_{i}(A\cup B)-
\sum\limits_{i=3}^{\Delta}i\sigma_{i}^{i}(A\cup B)-
\sigma_{2}^{2}(A\cup B)-\sigma^{2}_{2}(A)-\sigma^{2}_{2}(B)\\
&\ge -
\frac{|A\cup B|(|A\cup B|-1)}{2}-
\frac{|A|(|A|-1)}{2}-
\frac{| B|(|B|-1)}{2}.\end{split}
\end{equation*}
Thus
\begin{equation*}
\begin{split}
|B_{t}|& \geq\sum\limits_{f\in A\cup B}|f|-
\frac{|A\cup B|(|A\cup B|-1)}{2}-
\frac{|A|(|A|-1)}{2}-
\frac{| B|(|B|-1)}{2}\\
& \geq(|A|+|B|)s(H)-
\frac{(|A|+|B|)(|A|+|B|-1)}{2}-
\frac{|A|(|A|-1)}{2}-
\frac{| B|(|B|-1)}{2}\\
& =((1+s)\alpha-\alpha^{2}+(1+s)\beta-\beta^{2}-\alpha\beta)s(H)^{2}.
\end{split}
\end{equation*}
Let $f(\alpha,\beta)=
(1+s)\alpha-\alpha^{2}+(1+s)\beta-\beta^{2}-\alpha\beta$. By the same argument as above, we have
$$f(\alpha,\beta)\geq\left\{
\begin{array}{ll}
s+\frac{1}{4}  &  \mbox{when} \ s(H) \mbox{ is even,}\\
\frac{-s^{2}}{4}+s+\frac{1}{4}  &  \mbox{when} \ s(H) \mbox{ is odd}.
\end{array}\right.$$
Thus
$$tw([H]_2)\geq\left\{
\begin{array}{ll}
\frac{1}{4}s(H)^{2}+s(H)-1  &  \mbox{when} \ s(H) \mbox{ is even,}\\
\frac{1}{4}s(H)^{2}+s(H)-\frac{5}{4}  &  \mbox{when} \ s(H) \mbox{ is odd},
\end{array}\right.$$
and we have the conclusion. \end{proof}

\vskip.2cm

When $s(H)$ is even, let $H=(C_n^k)^*$, where $C_n^k$ is the $k^{th}$-power of an $n$-vertex cycle with vertex set $\{v_1,v_2,\ldots,v_n\}$ and edge set $\{(v_i,v_j)~|~\min\{|i-j|,i+n-j\}\leq k, 1\leq i<j\leq n\}$.  We can see that in this case $H$ is a 2-regular linear hypergraph and $s(H)=\delta(C_n^k)=2k$. By Theorem 1.3, $tw([H]_2)\geq \frac{1}{4}s(H)^{2}+s(H)-1=k^2+2k-1$. In \cite{Har}, Harvey and Wood showed that $tw(L(C_n^k))\leq k^2+2k-1$. Hence Theorem 1.3 is precisely sharp when $s(H)$ is even.

When $s(H)$ is odd, choose two matchings $X_1=\{1(n-k+1),2(n-k+2),\ldots,kn\}$ and $X_2=\{(k+1)(k+2),(k+3)(k+4),\ldots,(n-k-1)(n-k)\}$ in $C_n^k$. If $n$ is odd (resp. even), let $H=(C_n^k\setminus X_1)^*$ (resp. $H=(C_n^k\setminus(X_1\cup X_2))^*$). Then we have $s(H)=2k-1$. By Theorem 1.3, $tw([H]_2)\geq k^2+k-1$. In \cite{Har}, Harvey and Wood showed that
$$tw([H]_2)=tw(L(H^*))\leq\left\{
\begin{array}{ll}
k^2+k-1  &  \mbox{when} \ n \mbox{ is even,}\\
k^2+k-2  &  \mbox{when} \ n \mbox{ is odd}.
\end{array}\right.$$
Thus when $s(H)$ is odd, Theorem 1.3 is precisely sharp when $n$ is even; and within `+1' when $n$ is odd.

\section{ Upper bound}

\begin{proof}[of Theorem 1.4]
Let $(T,(B_{t})_{t\in T},(\lambda_{t})_{t\in T})$ be a supertree decomposition of a linear hypergraph $H$ with width $k$ such that $T$ has maximum degree at most 3. By the discussion in Section 1, we may assume that $r(H)\geq k-1$. (The existence of such a supertree decomposition $(T,(B_{t})_{t\in T},(\lambda_{t})_{t\in T})$ is well known and follows by a similar argument to Lemma 2.2.)

Say a hyperedge $f\in F$ is small if $|f|\leq k-1$ and large otherwise. For each $f\in F$, we  use $T(f)$ to denote the subtree of $T$ induced by $\lambda^{-1}(f)$. For each edge $e$ in $T$, let $A(e),B(e)$ denote the two  subtrees of $T- e$. If $e$ is also an edge of  $T(f)$ for some $f\in F(H)$, then let $A(e,f),B(e,f)$ denote two  subtrees of $T(f)- e$, where $A(e,f)\subseteq A(e)$ and $B(e,f)\subseteq B(e)$. For $t\in \lambda^{-1}(f)$, let $\gamma_{t}(f)=\{v\in V(H)|v\in f\cap g, g\in \lambda_t\setminus\{f\}\}$. Since $H$ is a linear hypergraph, we have $|\gamma_{t}(f)|\le k-1$. Denote  $\alpha(e,f)=\cup_{t\in A(e,f)}\gamma_{t}(f)$ and $\beta(e,f)=\cup_{t\in B(e,f)}\gamma_{t}(f)$. We have the following claim.
\vskip.2cm
{\noindent\bf Claim 2.}
For every large $f\in F$ there is an edge $e$ in $T(f)$ such that $|\alpha(e,f)|,|\beta(e,f)|\leq \frac{2}{3}|f|+\frac{1}{3}(k-1)$. 

\begin{proof}[of Claim 2] Assume for the sake of a contradiction that no such $e$ exists. Hence for all $e$ in $T(f)$, either $|\alpha(e,f)|$ or $|\beta(e,f)|$ is too ``large". Direct the edge $e$ towards $A(e,f)$ or $B(e,f)$ respectively. (If both $ |\alpha(e,f)|$, $|\beta(e,f)|$ are too large, then direct $e$ arbitrarily.) Given this orientation of $T(f)$, there must be a sink (all the edges incident to it direct to it), which we label $t_0$.

Let $e_{1},\ldots,e_{d}$ be the edges in $T(f)$ incident to $t_0$, where $d\in\{1,2,3\}$. Without loss of generality say that $e_{i}$ was directed towards $B(e_{i},f)$ for all $e_{i}$. Hence $|\beta(e_{i},f)|>\frac{2}{3}|f|+\frac{1}{3}(k-1)$ for all $i$.  Let $\alpha'(e_{i},t_0,f)=\alpha(e_{i},f)\setminus \gamma_{t_0}(f)$. Then $\alpha'(e_{i},t_0,f)\cap \alpha'(e_{j},t_0,f)=\emptyset$ when $i\neq j$.

 If $d=3$, then $\sum_{i=1}^{3}|\beta(e_{i},f)|>2|f|+k-1$. But $\sum_{i=1}^{3}|\beta(e_{i},f)|=
\sum_{i=1}^{3}\sum_{j\neq i}(|\gamma_{t_0}(f)|+|\alpha'(e_{j},t_0,f)|)=
2(|\gamma_{t_0}(f)|+\sum_{j=1}^{3}|\alpha'(e_{j},t_0,f)|)+|\gamma_{t_0}(f)|\leq
2|f|+k-1$, a contradiction.
If $d=2$, then $\sum_{i=1}^{2}|\beta(e_{i},f)|>\frac{4}{3}|f|+\frac{2}{3}(k-1)$. However, $\sum_{i=1}^{2}|\beta(e_{i},f)|=
\sum_{i=1}^{2}(\sum_{j\neq i}(|\alpha'(e_{j},t_0,f)|)+|\gamma_{t_0}(f)|)=
|\gamma_{t_0}(f)|+\sum_{j=1}^{2}|\alpha'(e_{j},t_0,f)|+|\gamma_{t_0}(f)|\leq
|f|+k-1$, a contradiction with $|f|>k-1$.
 If $d=1$, then $|\beta(e_{1},f)|>\frac{2}{3}|f|+\frac{1}{3}(k-1)$. However, $|\beta(e_{1},f)|=|\gamma_{t_0}(f)|\le k-1$, a contradiction with $|f|>k-1$.\end{proof}


\vskip.2cm
For each small hyperedge $f$ of $H$, arbitrarily select a base node in $\lambda^{-1}(f)$. For each large hyperedge $f$ of $H$, select an edge $e$ in the subtree $T(f)$ as guaranteed by Claim 2. Subdivide $e$ and declare the new node to be $b(f)$, the base node of $f$. If $e$ is selected for several different hyperedges, then subdivide it multiple times and assign a different base node for each hyperedge of $H$ that selected $e$. Denote the tree  after all of these subdivisions as $T'$. Together, this underlying tree $T'$ and the assignment $b$ gives a tree decomposition of $[H]_2$ in the same form as Lemma 2.1. Label the set of bags for this tree decomposition by $B'$. So the tree decomposition of $[H]_2$ is $(T',(B'_{t'})_{t'\in T'})$ and for each vertex $v\in V$, $B'^{-1}(v)=V(ST_v)$, where $ST_v$ is the subtree of $T'$ induced by $\cup_{f_i,f_j\in F(v)}Path(b(f_{i}),b(f_{j}))$. It remains to bound the width of this tree decomposition of $[H]_2$.

For each bag $B'_{t'}$ of $T'$, define a {\em corresponding bag} in $T$ as follows. If $t'\in T'$ is also in $T$, then the corresponding bag of $B'_{t'}$ is simply $\lambda_{t}$. If $t'\in T'$ is a subdivision node created by subdividing the edge $t_1t_2$, then the corresponding bag of $B'_{t'}$ is $\lambda_{t_1}$ or $\lambda_{t_2}$, chosen arbitrarily.

The following two claims give enough information to bound the width of $(T',(B'_{t'})_{t'\in T'})$.

\vskip.2cm
{\noindent\bf Claim 3.}
If $B'_{t'}$ is a bag of $T'$ with corresponding bag $\lambda_{t}$ ($t\in T$) and $v\in B'_{t'}$, then there is an edge $f\in F(v)$ such that $f\in \lambda_{t}$.

\begin{proof}[of Claim 3] Suppose that  $f\notin \lambda_{t}$ for all $f\in F(v)$. Then $t\notin \cup_{f\in F(v)} \lambda^{-1}(f)$.  If there are $f_i,f_j\in F(v)$ such that $T(f_{i})$ and $T(f_{j})$ are contained in different components of $T- t$, then $\lambda^{-1}(f_{i})\cap \lambda^{-1}(f_{j})=\emptyset$ a contradiction with $v\in f_{i}\cap f_{j}$. Thus for all $f\in F(v)$, $T(f)$ are all contained in the same component of $T- t$. Note that $b(f)$ is assigned inside of $\lambda^{-1}(f)$ (perhaps after some edges are subdivided, but this doesn't alter their positions relative to $\lambda_{t}$). Hence the subtree $ST_v$ in $T'$ doesn't include $t'$ which implies $v\notin B'_{t'}$, a contradiction. \end{proof}

\vskip.2cm

{\noindent\bf Claim 4.}
Suppose $f$ is a large hyperedge and $t'\in T'$ is not $b(f)$. If $\lambda_{t}$ ($t\in T$) is the corresponding bag of $B'_{t'}$,
then we have $|\gamma_{t}(f)|\leq\frac{2}{3}|f|+\frac{1}{3}(k-1)$.

\begin{proof}[of Claim 4]
Since $t'$ is not $b(f)$, there exists a component of $T'- b(f)$, say $T''$, containing $t'$. Let $v\in B'_{t'}\cap f$. Then $B'_{t'}$ is a bag in the subtree $ST_v$ in $T'$. Hence  there is $f'\in F(v)\setminus \{f\}$ such that $b(f')\in V(T'')$ since $v\in B'_{t'}$ .

Since $f$ is a large hyperedge, $b(f)$ is a subdivision node. Let $e$ be the edge in $T(f)$ that was subdivided to create $b(f)$. (The edge $e$ is also guaranteed by Claim 2.) Hence $V(T'')$ has non-empty intersection with exactly one of $V(A(e))$ and $V(B(e))$, say $V(T'')\cap V(A(e))\neq \emptyset$. Since $b(f')\in V(T'')$,  there must be $v\in \alpha(e,f)$  by $v\in f\cap f'$. Then $|\alpha(e,f)|\leq \frac{2}{3}|f|+\frac{1}{3}(k-1)$ by Claim 2. Hence  $B'_{t'}$ contains at most $\frac{2}{3}|f|+\frac{1}{3}(k-1)$ vertices in $f$, which means $\gamma_{t}(f)\leq\frac{2}{3}|f|+\frac{1}{3}(k-1)$. \end{proof}

\vskip.2cm
We now determine an upper bound on the size of a bag $B'_{t'},t'\in T'$. We count the vertices of $B'_{t'}$ by considering the number of vertices in a given hyperedge $f$ of $H$ contributes to $B'_{t'}$. By Claim 3, at most $k$ hyperedges of the corresponding bag $\lambda_{t}$ contribute  to $B'_{t'}$.

If $f$ is small, it contributes at most $k-1$ vertices to $B'_{t'}$.
If $f$ is large and $t'\neq b(f)$, by Claim 4, $f$ contributes at most $\frac{2}{3}r(H)+\frac{1}{3}(k-1)\geq k-1$ vertices to $B'_{t'}$.
If $f$ is large and $t'=b(f)$, then $f$ contributes at most $r(H)\geq \frac{2}{3}r(H)+\frac{1}{3}(k-1)$ vertices. Therefore, we conclude the highest possible contribution of a large $f$ with $t'=b(f)$ is greater than that of a large $f$ with $t'\neq b(f)$ which is greater than that of a small $f$.
Note that $t'=b(f)$ for at most one $f\subseteq F$. Hence
$$|B'_{t'}|\leq (k-1)\left(\frac{2}{3}r(H)+\frac{1}{3}(k-1)\right)+r(H)=
\frac{2}{3}kr(H)+\frac{1}{3}(k-1)^{2}+\frac{1}{3}r(H),$$
where there exists $k-1$

If we set $(T,(B_{t})_{t\in T},(\lambda_{t})_{t\in T})$ to be a minimum width hypertree decomposition of $H$, then $k=stw(H)$ and so
$$tw([H]_2)\leq\frac{2}{3}stw(H)r(H)+\frac{1}{3}(stw(H)-1)^{2}+\frac{1}{3}r(H)-1.$$ \end{proof}

{\noindent\bf Remark} If $tw(L(H))\leq c$ for some constant $c$, we can get the exact value of $stw(H)$ in polynomial time by  $stw(H)=tw(L(H))+1$. Then Theorem 1.4 can give a useful upper bound for $tw([H]_2)$.

\section{ Conclusion}

Treewidth is an important graph parameter in structural graph theory  and in algorithmic graph theory. This paper studies the treewidth of corresponding graphs of linear hypergraphs. Let $G=(V,E)$ be a graph. There is a linear hypergraph $H=(V,F)$ such that $[H]_2\cong G$. For example, we first find cliques $C_1,C_2,\ldots,C_s$ in $G$ such that $\cup_{i=1}^sV(C_i)=V$ and $|V(C_i)\cap V(C_j)|\le 1$ for $i\not=j$ and $1\le i,j\le s$. Now we construct a hypergraph $H=(V,F)$ with $V(H)=V(G)$ and $F=\{C_1,\ldots,C_s\}\cup \{\{x,y\}|x,y\in V,xy\in E(G)\setminus \cup_{i=1}^sE(C_i)\}$. Obviously, $H$ is a linear hypergraph and $[H]_2\cong G$. Thus we provide a method to estimate the  bound of treewidth of graph by  the parameters of the hypergraph. 

\acknowledgements
This work is partially supported by  the National Natural Science Foundation of China (Grant 11771247 \& 11971158) and  Tsinghua University Initiative Scientific Research Program. Many thanks to the anonymous referees for their many helpful comments and suggestions, which have considerably improved the presentation of the paper.

\bigskip
{\bf Appendix A}

We calculate the minimum value of $f(\alpha,\beta)=\frac{\alpha}{\Delta}+\frac{\beta}{\Delta}+\frac{1}{\delta(\delta-1)}(-\alpha^{2}+s\alpha-\beta^2+s\beta)
$, under the condition $s\leq\alpha,\beta\leq\frac{1}{2}$, $\alpha+\beta>\frac{1}{2}$, $0<s\leq\frac{1}{2}$ and $\Delta\geq \delta\geq2$.

Note that $f(\alpha,\beta)=f(\beta,\alpha)$. Consider the partial derivative of $f(\alpha,\beta)$ with respect to $\alpha$ and $\beta$, we have
\[\frac{\partial f(\alpha,\beta)}{\partial \alpha} = \frac{1}{\Delta}-\frac{2\alpha}{\delta(\delta-1)}+\frac{s}{\delta(\delta-1)}, ~~\frac{\partial f(\alpha,\beta)}{\partial \beta} = \frac{1}{\Delta}-\frac{2\beta}{\delta(\delta-1)}+\frac{s}{\delta(\delta-1)}.\]
In the boundary $\alpha+\beta=\frac{1}{2}$, we have
$$f\left(\alpha,\frac{1}{2}-\alpha\right)=
\frac{-8\Delta\alpha^{2}+4\Delta\alpha+2\delta^{2}-2\delta-\Delta+2\Delta s}{4\Delta\delta(\delta-1)}$$
and
$$\frac{\partial f\left(\alpha,\frac{1}{2}-\alpha\right)}{\partial \alpha} =\frac{1-4\alpha}{\delta(\delta-1)}.$$Thus we have $f(\alpha,\beta)\ge \min\left\{f\left(\frac{1}{2},\frac{1}{2}\right),f\left(\frac{1}{2},s\right),f\left(\frac{1}{2}-s,s\right)\right\}$.
Note that
\begin{eqnarray*}
f\left(\frac{1}{2},\frac{1}{2}\right)&=&\frac{2\Delta s+2\delta^{2}-2\delta-\Delta}{2\Delta(\delta-1)},\\
f\left(\frac{1}{2},s\right)&=&\frac{2\Delta s-2\delta-\Delta-4\delta s+4\delta^{2}s+2\delta^{2}}{4\Delta(\delta-1)},\\
f\left(\frac{1}{2}-s,s\right)&=&\frac{-8\Delta s^{2}+6\Delta s+2\delta^{2}-2\delta-\Delta}{4\Delta(\delta-1)}.
\end{eqnarray*}
Since  $\Delta\le 2\delta^{2}-2\delta$, $f(\frac{1}{2},\frac{1}{2})\ge f(\frac{1}{2},s)$. Thus $f(\alpha,\beta)\ge \min\left\{f\left(\frac{1}{2},s\right),f\left(\frac{1}{2}-s,s\right)\right\}$.

{\bf Appendix B}

Here we show that
$-\sum\limits_{j=2}^{\Delta}\sum\limits_{i=j+1}^{\Delta}\sigma^{j}_{i}(A\cup B)(j-1)-
\sum\limits_{i=\delta}^{\Delta}i\sigma_{i}^{i}(A\cup B)\ge -
\frac{|A\cup B|(|A\cup B|-1)}{2}$.

We consider the following  linear programming
\begin{equation*}
\begin{split}
 {\bf min}\quad &-\sum\limits_{j=2}^{\Delta}\sum\limits_{i=j+1}^{\Delta}(j-1)\sigma_{i}^j(A\cup B)-
\sum\limits_{i=\delta}^{\Delta}i\sigma_{i}^i(A\cup B)\\
 {\bf s.t}\quad &\sum\limits_{j=2}^{\Delta}\sum\limits_{i=j}^{\Delta}\frac{j(j-1)}{2}\sigma_{i}^j(A\cup B)\leq \frac{n'(n'-1)}{2},\\
 &-\sigma_{i}^j(A\cup B)\leq 0,\quad 2\leq j\leq i\leq \Delta,
\end{split}
\end{equation*}
where $n'=\left|A\cup B\right|$. The $KKT$ condition of this linear programming is
\begin{equation*}
\begin{split}
& 1-j+u\frac{j(j-1)}{2}-u_{i,j}=0,\quad 2\leq j<i\leq \Delta,\\
& -i+u\frac{i(i-1)}{2}-u_{i,i}=0,\quad i=\delta,\ldots,\Delta,\\
& u,u_{i,j}\geq 0,\quad 2\leq j\leq i\leq \Delta,\\
& u\left(\sum\limits_{j=2}^{\Delta}\sum\limits_{i=j}^{\Delta}\frac{j(j-1)}{2}\sigma_{i}^j(A\cup B)- \frac{n'(n'-1)}{2}\right)=0,\\
& u_{i,j}\sigma_{i}^j(A\cup B)=0,\quad 2\leq j\leq i\leq \Delta.
\end{split}
\end{equation*}
We can easily get $u=1$, $u_{i_0,2}=0$ and $\sigma_{i_{0}}^2(A\cup B)=\frac{n'(n'-1)}{2}$ for some $\delta\leq i_{0}\leq \Delta$, and any other $u_{i,j}=\sigma_{i}^j(A\cup B)=0$.
Thus $-\sum\limits_{j=2}^{\Delta}\sum\limits_{i=j+1}^{\Delta}\sigma^{j}_{i}(A\cup B)(j-1)-
\sum\limits_{i=\delta}^{\Delta}i\sigma_{i}^{i}(A\cup B)\ge -
\frac{|A\cup B|(|A\cup B|-1)}{2}$.

{\bf Appendix C}

Here we show the lower bound of
$$-\sum\limits_{j=2}^{\Delta}\sum\limits_{i=j+1}^{\Delta}(j-1)\sigma^{j}_{i}(A\cup B)-
\sum\limits_{i=3}^{\Delta}i\sigma_{i}^{i}(A\cup B)-
\sigma_{2}^{2}(A\cup B)-\sigma^{2}_{2}(A)-\sigma^{2}_{2}(B).$$

We consider the following  linear programming
\begin{equation*}
\begin{split}
 {\bf min}\quad &-\sum\limits_{j=2}^{\Delta}\sum\limits_{i=j+1}^{\Delta}(j-1)\sigma_{i}^j(A\cup B)-
\sum\limits_{i=\delta}^{\Delta}i\sigma_{i}^i(A\cup B)-
\sigma_{2}^2(A\cup B)-\sigma^{2}_{2}(A)-\sigma^{2}_{2}(B)\\
 {\bf s.t}\quad &\sum\limits_{j=2}^{\Delta}\sum\limits_{i=j}^{\Delta}\frac{j(j-1)}{2}\sigma_{i}^j(A\cup B)\leq \frac{n'(n'-1)}{2},\\
 &\sigma^{2}_{2}(A)+\sigma^{2}_{2}(B)\leq \sigma_{2}^2(A\cup B),\\
 &0\leq \sigma^{2}_{2}(A)\leq \frac{n'_{1}(n'_{1}-1)}{2},\\
&0\leq \sigma^{2}_{2}(B)\leq \frac{n'_{2}(n'_{2}-1)}{2},\\
 &-\sigma_{i}^j(A\cup B)\leq 0,\quad 2\leq j\leq i\leq \Delta,
\end{split}
\end{equation*}
where $n'=\left|A\cup B\right|$, $n'_{1}=|A|$ and $n'_{2}=|B|$. The $KKT$ condition of this linear programming is
\begin{equation*}
\begin{split}
& 1-j+u_{1}\frac{j(j-1)}{2}-u_{i,j}=0,\quad 2\leq j<i\leq \Delta,\\
& -i+u_{1}\frac{i(i-1)}{2}-u_{i,i}=0,\quad i=3,\ldots,\Delta,\\
&-1+u_{1}-u_{2}-u_{2,2}=0, \\
&-1+u_{2}+u_{3}-u_{5}=0, \\
&-1+u_{2}+u_{4}-u_{6}=0,\\
& u_{i,j}\geq 0,\quad 2\leq j\leq i\leq \Delta,\\
& u_{l}\geq 0,\quad l=1,\ldots,5,\\
& u_{1}(\sum\limits_{j=2}^{\Delta}\sum\limits_{i=j}^{\Delta}\frac{j(j-1)}{2}\sigma_{i}^j(A\cup B)- \frac{n'(n'-1)}{2})=0,\\
& u_{2}(\sigma_2^2(A)+\sigma_2^2(B)-\sigma_{2}^2(A\cup B))=0,~
u_{3}(\sigma_2^2(A)-\frac{n'_{1}(n'_{1}-1)}{2})=0,\\
& u_{4}(\sigma_2^2(B)-\frac{n'_{2}(n'_{2}-1)}{2})=0,
u_{5}\sigma_2^2(A)=0,
u_{6}\sigma_2^2(B)=0,\\
& u_{i,j}\sigma_{i}^j(A\cup B)=0,\quad 2\leq j\leq i\leq \Delta.
\end{split}
\end{equation*}
We can get $\sigma_{2}^2(A\cup B)=\frac{n'(n'-1)}{2}$, $\sigma_2^2(A)=\frac{n'_{1}(n'_{1}-1)}{2}$, $\sigma_2^2(B)=\frac{n'_{2}(n'_{2}-1)}{2}$ and all other $\sigma_{i}^j(A\cup B)=0$. Thus
\begin{equation*}
\begin{split}&-\sum\limits_{j=2}^{\Delta}\sum\limits_{i=j+1}^{\Delta}(j-1)\sigma^{j}_{i}(A\cup B)-
\sum\limits_{i=3}^{\Delta}i\sigma_{i}^{i}(A\cup B)-
\sigma_{2}^{2}(A\cup B)-\sigma^{2}_{2}(A)-\sigma^{2}_{2}(B)\\
&\ge -
\frac{|A\cup B|(|A\cup B|-1)}{2}-
\frac{|A|(|A|-1)}{2}-
\frac{| B|(|B|-1)}{2}.\end{split}
\end{equation*}

\end{document}